\documentclass[12pt]{amsart}

\usepackage{amsmath, amssymb, amsthm, amscd, amsfonts}

\newcommand{\Sk}{S_{2k}(\GN)}
\newcommand{\Mk}{M_{2k}(\GN)}
\newcommand{\Minf}{M_{2k}^{\infty}(\GN)}
\newcommand{\SL}{\text{SL}_{2}(\mathbb{Z})}
\newcommand{\SLR}{\text{SL}_{2}(\mathbb{R})}
\newcommand{\GN}{\Gamma}
\newcommand{\JN}{j_{\GN}}
\newcommand{\DN}{\Upsilon}
\newcommand{\F}{\mathcal{F}}
\newcommand{\bF}{\partial\F}
\newcommand{\XN}{\GN\setminus\UH^*}
\newcommand{\Ek}{E_{2k}^{\GN}}
\newcommand{\UH}{\mathfrak{H}}
\newcommand{\w}{\omega}

\numberwithin{equation}{section}

\newtheorem{theorem}{Theorem}[section]
\newtheorem{lemma}[theorem]{Lemma}
\newtheorem{corollary}[theorem]{Corollary}
\newtheorem*{definition}{Definition}
\newtheorem{proposition}[theorem]{Proposition}
\theoremstyle{remark}

\newtheorem*{remark}{Remark}

\begin{document}

\title{On zeros of Eisenstein series for genus zero Fuchsian groups}
\author{Heekyoung Hahn}
\address{Department of Mathematics, University of Rochester, Rochester, NY 14627 USA}
\email{hahn@math.rochester.edu}

\begin{abstract}
Let $\GN\leq\SLR$ be a genus zero Fuchsian group of the first kind with $\infty$ as a cusp, and let $\Ek$ be the holomorphic Eisenstein series of weight $2k$ on $\GN$ that is nonvanishing at $\infty$ and vanishes at all the other cusps (provided that such an Eisenstein series exists). Under certain assumptions on $\GN,$ and on a choice of a fundamental domain $\F$, we prove that all but possibly $c(\GN,\F)$ of the non-trivial zeros of $\Ek$ lie on a certain subset of $\{z\in\mathfrak{H}\,:\,\JN(z)\in\mathbb{R}\}$. Here $c(\GN,\F)$ is a constant that does not depend on the weight $2k$ and $\JN$ is the canonical hauptmodul for $\GN.$ 
\end{abstract}

\thanks{Research is supported in part by a National Science Foundation FRG grant (DMS 0244660)}
\maketitle

\section{Introduction and statement of the result.}

Recently, Rudnick \cite{Rudnick} showed that, under the Generalized Riemann Hypothesis, the zeros of Hecke eigenforms for the full modular group $\SL$ become equidistributed with respect to the hyperbolic measure on a fundamental domain for the action of $\SL.$ However, the situation is quite different in case of zeros of Eisenstein series for $\SL.$ Rankin and Swinnerton-Dyer \cite{Rankin-Swinnerton} proved that all the zeros of Eisenstein series of weight $k\geq 4$ for $\SL$ in the standard fundamental domain lie on $A:=\{ e^{i\theta}| \frac{\pi}{2}\leq \theta \leq \frac{2\pi}{3}\}.$ Rankin \cite{Rankin-Poincare} proved that similar properties hold for certain Poincar\'e series on $\SL.$ Getz \cite{Getz} generalized Rankin and Swinnerton-Dyer's result by providing conditions under which the zeros of other modular forms lie only on the arc $A.$ Given these work, it is natural to investigate the behavior of the zeros of families of modular forms (in particular, the family of even weight Eisenstein series) on other groups. The goal of the present paper is to find the location of zeros of certain Eisenstein series associated with genus zero Fuchsian group.

 From the outset we denote by $\GN$ a Fuchsian group of the first kind \cite{Shimura}, by $\UH$ the upper half-plane, and by $\UH^*$ the union of $\UH$ with the set of cusps of $\GN.$ We say that a meromorphic function $f$ on $\UH$ is a {\it meromorphic modular form} of {\it weight k} for $\GN$ if
\begin{equation}
f\Big(\frac{az+b}{cz+d}\Big)=(cz+d)^{-k}f(z)
\end{equation}for all $z\in\UH$ and $\left(\begin{smallmatrix} a & b \\ c & d\end{smallmatrix}\right)\in \GN,$ and $f$ is meromorphic at the cusps. If $k=0,$ then $f$ is known as a {\it modular function} on $\GN.$ Further, we say that $f$ is a {\it holomorphic modular form} if $f$ is holomorphic on $\UH$ and holomorphic at cusps. A holomorphic modular form is said to be a {\it cusp form} if it vanishes at the cusps of $\GN.$ As usual, we denote by $\Mk$ (resp. $\Sk$) the space of holomorphic modular forms (resp. cusp forms) of weight $2k$ for $\GN.$ Moreover denote by $\Minf$ the space of weakly holomorphic modular forms on $\GN$ (i.e. holomorphic on $\UH$ but not necessarily at the cusps). 
 
If $\GN$ has a cusp at infinity with width $h,$ then each $f\in\Mk$ has the Fourier expansion at infinity
\begin{equation}\label{q}
f(z)=\sum_{m=0}^{\infty}a_m q_{h}^{m}, \quad q_{h}:=e^{2\pi i z/h},\quad z\in\UH.
\end{equation}Note that a modular form can be {\it identified} with its $q$-expansion. we say that the form has {\it real Fourier coefficients} if $a_m$ is real for all $m$.

Assume that $\GN$ has $\infty$ as a cusp.  Let $\Ek \in \Mk$ be the unique modular form that is orthogonal to all cusp forms vanishes at all cusps but $\infty$ and has constant term $1$ in its Fourier expansion at $\infty$ (provided that such a modular form exists).  We will refer to $\Ek$ as the Eisenstein series associated to the cusp $\infty$.  This is the Eisenstein series to which we referred in the title of this paper.
 \begin{remark}
Note that such form doesn't always exist. For instance, take $k=1$ and $\GN=\SL.$ 
\end{remark}
\begin{definition}
The group $\GN$ is {\bf good} for the weight $2k$ if
\begin{enumerate}
\item{The Eisenstein series $\Ek$ has real Fourier coefficients (if it exists).}
\item{The space $\Mk$ has a basis of forms with real Fourier coefficients.}
\end{enumerate}
\end{definition}
\begin{remark} For example, $\GN_0(N)$ is good for all $N$.  See \S \ref{examples}.
\end{remark}

Assume $\GN$ has genus zero. In this case, it is well-known that there is a unique weakly holomorphic modular form$$\JN(z)=\frac{1}{q}+\sum_{n=1}^{\infty}a_{n}q^{n}\,\in\,M_0^{\infty}(\GN).$$ This function $\JN$ is referred as the {\it canonical hauptmodul}.

Suppose $\GN$ has width $h$ at the cusp infinity. We choose a fundamental domain $\F$ for $\GN$ contained in $\{z \,|\,-h/2\leq \Re(z)\leq h/2\}.$ Define $y_0\in \mathbb{R}$ by
\begin{equation}\label{y0}
y_{0}:=\inf\, \{ y\,|\,\pm h/2+iy \in \bF\},
\end{equation}where $\bF$ denotes the boundary of $\F.$ We refer to the complement of the vertical lines 
$$L:=\{z \,|\, \Re(z)=\pm h/2, \,\text{and}\, \Im(z)\geq y_0 \} \subset \bF$$ 
as the \emph{lower arcs $A$}. 

For the purpose of our paper, we need to place some restrictions on our fundamental domain. Assume that we can parametrize $A$ by a piecewise smooth function 
\begin{eqnarray*}
z_A=x_A+iy_A:(0,1)&\longrightarrow& A\\
t\,\, &\longmapsto &z_A(t).
\end{eqnarray*}
Fixing such parametrization, consider the set of points
\begin{equation}\label{Crit}
Crit(\F):=\Big\{t \in (0,1): \frac{y'(t)}{(\JN\circ z_A)'(t)}=0 \textrm{ or } \frac{y'(t)}{(\JN\circ z_A)'(t)}\textrm{ is undefined}\Big\}.
\end{equation}Here the $'$ indicates differentiation with respect to $t.$  We say that $\frac{y'(t)}{(\JN\circ z_A)'(t)}$ {\it changes sign} at $t_0$ if, for any sufficiently small $\epsilon_1,\epsilon_2 >0,$
$$\Big(\frac{y'(t_0-\epsilon_1)}{(\JN\circ z_A)'(t_0-\epsilon_1)}\Big)
\Big(\frac{y'(t_0+\epsilon_2)}{(\JN\circ z_A)'(t_0+\epsilon_2)}\Big)<0.$$We make the following definition:
\begin{definition}
A fundamental domain $\F$ for $\GN$ is {\bf acceptable} if we can parametrize $A$ by a piecewise smooth path $z_A$ as above and the following conditions hold:
\begin{enumerate}
\item{The hauptmodul $\JN$ is real on $\bF.$}
\item{The set $Crit(\F)$ is a discrete set.}
\item{If $\gamma z_A(t_0)=z_A(\hat{t}_0)$ for some $\gamma\in\GN$ and for some $t_0,\,\hat{t}_0 \in Crit(\F)$, then $\frac{y'(t)}{(\JN\circ z_A)'(t)}$ changes sign at $t_0$ if and only if it changes sign at $\hat{t}_0.$}
\end{enumerate}
\end{definition}
\begin{remark}
These appear to be natural restrictions on $\F$.   In \S \ref{examples} we deal with the case $\GN=\GN_0(3)$, for example.  We must point out, however, that we do not have a proof that such a fundamental domain exists for arbitrary good genus zero $\GN$.  
\end{remark}
Now, define 
\begin{equation}\label{constant}
c(\GN,\F):=\Big|\big\{ z(t_0): t_0 \in Crit(\mathcal{F}) \textrm{ and } \frac{y'(t)}{(\JN\circ z_A)'(t)} \textrm{ changes sign at }\,t_0\big\}/\sim\Big|,
\end{equation}where we consider two points in the set to be equivalent if they are $\GN$-conjugates of eachother. Note that  $c(\GN,\F)$ depends on $\GN$ and on a choice of $\F$, but not on weight $2k.$ 

For each $f \in M_{2k}(\GN)$, let 
\begin{equation}\label{Poly}
P(f,X):=\prod_{\substack{Q \in \mathcal{F}\setminus\{\infty\}\\\textrm{ord}_Q f\neq 0}} (X-\JN(Q))^{\textrm{ord}_Q(f)-e_Q\{k(1-1/e_Q)\}}.
\end{equation}
Here $\{k(1-1/e_Q)\}$ is a fractional part of $k(1-1/e_Q),$ and $e_Q:=|\mathrm{Stab}_{\overline{\GN}}(Q)|,$ where $\overline{\GN}:=\GN/\{\pm 1\}$ if $-1\in\GN,$ and $\overline{\GN}=\GN$ otherwise. Thus, if $P(f,\JN(Q))=0$, then $f(Q)=0,$ with the converse being true if $Q$ is not an elliptic fixed point. It follows from Lemma \ref{trivialzero} below that $P(f,X) \in \mathbb{C}[X]$ (this is obvious if the divisor of $f$ contains no elliptic fixed points).  We call $P(f,X)$ the \emph{divisor polynomial of $f$}.  We can now state our main theorem:
\begin{theorem}\label{main}
Let $\Gamma$ be a genus zero group that is good for the weight $2k\geq 2.$  Suppose that we can choose an acceptable fundamental domain $\F$ for $\GN$, and moreover that $\JN$ has real Fourier coefficients.  If $\Ek$ exists, then all but possibly $c(\GN, \F)$ of the zeros of $P(\Ek,X)$ are real, simple, and contained in the interval $[\JN(-\frac{h}{2}+iy_0),\,\infty),$ where $y_0$ is defined as \eqref{y0}. 
\end{theorem}

\begin{remark} In the case $\GN=\SL$ and $\mathcal{F}$ is the standard fundamental domain bounded by the lines at $\Re(z)= \pm \frac{1}{2}$ and $|z|=1$, this theorem implies that the standard Eisenstein series on $\mathrm{SL}_2(\mathbb{Z})$ has all of its zeros in $\mathcal{F}$ on the imaginary axis and the curve $|z|=1$.  Thus Theorem \ref{main} is weaker than Rankin and Swinnerton-Dyer's result \cite{Rankin-Swinnerton} in this case, which proves that the zeros all lie on arc $\mathcal{F} \cap \{z:|z|=1\}$, but it has the advantage that the same argument applies to a family of groups (see \S \ref{examples}, for example).
\end{remark}

Recalling the elementary fact that the complex zeros of a polynomial with real coefficients must occur in pairs, we obtain the following immediate corollary of Theorem \ref{main}.

\begin{corollary}
Assume the hypotheses of Theorem \ref{main}.  If $c(\GN,\F)$ is odd, then all but possibly $c(\GN,\F)-1$ of the zeros of $P(E_{2k}^{\GN}, X)$ are real.  
\end{corollary}

The proof of Theorem \ref{main} is given in \S \ref{proof1}.  It can be viewed as an adaptation of the classical argument proving that the zeros of an orthogonal polynomial are all real and simple (see \cite[Theorem 5.4.1]{AAR}, for example).  We now outline the argument.  Using the Riemann-Roch theorem, we first construct a modular form $\DN\in \Mk$ with a zero of order $\dim(\Mk)-1$ at infinity and no other ``nontrivial'' zeros. Suppose for a contradiction that  the result is not true. Then we can construct a polynomial $Q(X)$  of degree less than $\dim(\Mk)-1$ whose zeros are exactly the zeros of $\Ek$ that occur with odd multiplicity and the image under $j_{\Gamma}$ of points where $\frac{y'(t)}{(\JN\circ z_A)'(t)}$ changes sign along the lower arc $A.$  Then $Q(j_{\Gamma})\Upsilon$ will be an element of $\Mk$ that vanishes at $\infty$, and it follows that the Petersson inner product $\langle\Ek,Q(\JN)\DN\rangle$ is zero. On the other hand, using Green's theorem, we relate the Petersson inner product $\langle\Ek, Q(\JN)\DN\rangle$ to a line integral along the lower arcs. By our construction of $Q(X),$ we show that this line integral  cannot be zero, resulting in the desired contradiction.

We close this introduction with the following remark:

\begin{remark} Suppose that a genus zero $\Gamma$ is good and admits an acceptable fundamental domain.   Since $\Gamma$ is of genus zero, $j_{\Gamma}$ defines a model  $j_{\Gamma}:\Gamma \backslash \UH^*\to \mathbb{P}^1(\mathbb{C})$ of $\Gamma \backslash \UH^*$ as the Reimann sphere.  One can interpret Theorem \ref{main} as saying that all but a bounded number of the zeros of $\Ek$ lie on the real points of this model.  It would be interesting to see if the some modification of the ideas of this paper could prove that, for suitable higher genus congruence groups $\Gamma$, some proportion of the zeros of corresponding Eisenstein series lie either at the Weierstrass points or the real points of a suitable model for $\Gamma \backslash \UH^*$.  We should emphasize that this is only a guess as to what occurs in these cases; we have not even done any numerical computations.  
\end{remark}

\section{Examples} \label{examples}

In this section, we give some examples of good groups and acceptable fundamental domains. We note that we found the fundamental domain drawing applet (see \cite{Verrill}) useful in producing these examples.  

Clearly $\SL$ is good. The {\it principal congruence subgroup} $\GN(N)$ of $\SL$ 
$$\GN(N):=\Big\{\gamma \in \SL\,:\, \gamma \equiv \begin{pmatrix} 1 & 0 \\ 0 & 1\end{pmatrix} \pmod{N}\Big\}$$ 
is good for all $k$ \cite[p.163]{Schoeneberg}. For integer $N\geq 1,$ any group $\GN$ satisfying 
$\GN_{1}(N)\leq \GN \leq\GN_{0}(N)$ is good \cite{Diamond-Im} for all $k,$ where 
$$\GN_{1}(N):=\Big\{\gamma \in \SL\,:\, \gamma \equiv \begin{pmatrix} 1 & * \\ 0 & 1\end{pmatrix} \pmod{N}\Big\},$$
$$\GN_{0}(N):=\Big\{\gamma \in \SL\,:\, \gamma \equiv \begin{pmatrix} * & * \\ 0 & *\end{pmatrix} \pmod{N}\Big\}.$$
Define the group $\GN_0^*(N)\leq\SLR$ as in \cite[p.27]{Shimura} to be the extension of $\GN_{0}(N)$ by the Friche involution $$W_{N}:=\begin{pmatrix} 0 & -\sqrt{N}^{-1}\\ \sqrt{N} & 0\end{pmatrix}.$$ Since $\Gamma_0(N)$ is good for all $2k$ and an involution has eigenvalues $\pm 1$, it follows that there is a basis of eigenvectors for $W_N$ acting on $M_{2k}(\Gamma_0(N))$ whose Fourier coefficients are all real.  The forms in this basis with eigenvalue $+1$ will be a basis for $M_{2k}(\Gamma_0^*(N))$; it follows that the  $\GN_0^*(N)$ are good.  

Some good groups have the {\it acceptable} fundamental domains.  Before we consider any examples, we point out that if $j_{\Gamma}$ has real Fourier coefficients, then it is \emph{always} real on $\Re(z)=\pm \frac{h}{2}$, where $h$ is the width of the cusp at $\infty$.  This is easy to see from its Fourier expansion at $\infty$. 

We now consider a few examples.  If $\Gamma=\SL$, the classical fundamental domain bounded by $\Re(z)= \pm \frac{1}{2}$ and the unit circle is acceptable (see \cite[pp. 31--32]{Schoeneberg}), and $c(\Gamma,\F)=0$.

Consider $\GN=\GN_0(3).$ In this case, for $k \geq 2$ the Eisenstein series $\Ek$
has Fourier expansion  
\begin{equation}\label{Ek-3}
\Ek(z)=1-\frac{8k}{B_{2k}(3^{2k}-1)}\sum_{n=1}^{\infty}\delta_{2k-1}(n)q^n,
\end{equation}where $B_{2k}$ is the $2k$th Bernoulli number, 
 $$
 \delta_{2k-1}(n):=\sum_{d|n}\alpha(d)d^{2k-1},
 $$ 
 and
\begin{equation}
\alpha(d):=\begin{cases} 1, &\text{if }\, 3\equiv 0 \pmod{d}\\
 -1/2, &\text{if }\, 3\not\equiv  0 \pmod{d}.\end{cases}
\end{equation} 
(see \cite[Theorem 1.1]{Boylan}).
Denote by $\F_{3}$ the fundamental domain bounded by vertical lines $\{z\,:\,\Re(z)=\pm \frac{1}{2}, \Im(z)\geq \frac{\sqrt{3}}{6}\}$ and lower arcs $A_{l}:=\{z: -\frac{1}{2}\leq\Re(z)\leq 0, \text{ and }|z+\frac{1}{3}|=\frac{1}{3}\}$ and $A_{r}:=\{z: 0\leq \Re(z)\leq\frac{1}{2}, \text{ and }|z-\frac{1}{3}|=\frac{1}{3}\}.$ Let for $q:=e^{2\pi i z},$
$$
j_{3}(z):=q^{-1}+783q+8672q^{2}+65367q^{3}+371520q^{4}+\cdots
$$ 
be the canonical Hauptmodul for  $\GN_{0}(3).$ It is clear from the Fourier expansion at $\infty$ that $j_{3}(z)=\overline{j_{3}(-\bar{z})}.$  For $z \in A_{l},$ a simple calculation shows that 
$\gamma z=-\bar{z},$ where 
$$
\gamma:=\begin{pmatrix} 1 & 0 \\ -3 & 1\end{pmatrix}\in\GN_{0}(3),
$$
 and hence $j_{3}(z)=\overline{j_{3}(-\bar{z})}=\overline{j_{3}(\gamma z)}=\overline{j_{3}(z)},$ because $j_{3}$ is invariant under $\GN_{0}(3).$ For $z\in A_{r},$ we have $\gamma'z=z$ where
 $$
 \gamma':=\begin{pmatrix} 1 & 0 \\ 3 & 1\end{pmatrix}\in\GN_{0}(3).
 $$
 Therefore we have $j_{3}(z)\in\mathbb{R}$ for $z\in\bF_{3}.$ Since $c(\GN_0(3), \F_3)=1,$ all zeros of $P(\Ek,X)$ are real, and all but possibly  one zero are simple and lie on the interval $[-42,\infty)$.  Here we use the fact that  $j_{3}(\frac{-3+i\sqrt{3}}{6})=-42$ (see \cite[p. 317]{Ono-Bring}).

\section{Prelimiaries on divisor polynomials}\label{sec2}

In this section we prove a proposition which gives an alternative construction of the {\it divisor polynomial} $P(f,X)$ attached to a modular form $f\in\Mk.$ The idea is basically that of Ono (see \cite[Section 2.6.]{Ono}), but things are complicated in our case due to the fact that the algebra structure of $\bigoplus_{k\geq 0}\Mk$ is not always as simple as it is in the case $\GN=\SL.$ To overcome this difficulty, we rely on a few easy applications of Riemann-Roch theory. By abuse of language, we sometimes call a point of $\XN$ an {\it elliptic point} or a {\it cusp}, if it corresponds to an elliptic point or a cusp on $\UH^*$ with respect to $\GN.$
\begin{lemma}\label{trivialzero}
Any modular form in $\Mk$ has zeros of oder at least $e_Q\{k(1-1/e_Q)\}$ at every elliptic points $Q\in \UH^{*}.$ 
\end{lemma}
\begin{proof}
Let $f \in \Mk$ and let $\w$ be the corresponding $k$-fold differential form on $\XN.$ Then if $Q$ is an elliptic fixed point, then
\begin{equation}\label{P-Q}
\text{ord}_{Q}(f)=e_Q\text{ord}_{P}(\w)+k(e_Q-1).
\end{equation}Here $P=\pi(Q),$ where $\pi : \UH^* \rightarrow \XN$ is the canonical projection (See \cite[Lemma 4.11]{Milne}). Since $Q$ is elliptic, it is well-known that $e_Q=2$ or $3.$ In view of this, it is easy to check that \eqref{P-Q} implies the lemma.
\end{proof}
Let
\begin{equation}\label{d}
d:=\dim(\Mk)
\end{equation}for a good $\GN$ of fixed weight $k\geq 1$ with genus zero. If $e_Q$ is defined as \eqref{Poly} for every elliptic point $Q$ on $\UH^*$ with respect to $\GN,$ then we have the following proposition.
\begin{proposition}\label{Delta}
If $\Gamma$ is of genus zero, then there exists a modular form $\DN\in \Mk$ with a zero of order $d-1$ at infinity, zeros of order $e_Q\{k(1-1/e_Q)\}$ at the elliptic point $Q$ and no other zeros. Moreover, if $\Gamma$ is good,  then all the Fourier coefficients of $\DN$ are real.
\end{proposition}

We defer the proof of Proposition \ref{Delta} to the end of the section.  Assuming it for the moment, we claim that if $f \in \Mk$ has $1$ as the leading term of its Fourier expansion at $\infty$, then 
\begin{eqnarray}
f \DN^{-1}=P(f, \JN)
\end{eqnarray}
where $P(f,X)$ is the divisor polynomial of $f$, defined as in \eqref{Poly}.  Indeed, by Proposition 
\ref{Delta}, $f \DN^{-1} \in M_0^{\infty}(\GN)$ and has the zeros specified by \eqref{Poly} to the correct order.  Since every element of $M_0^{\infty}(\GN)$ is a polynomial in $\JN$, the claim follows.  

In the following proof, we use standard notation from Riemann-Roch theory on algebraic curves.  See (\cite[pp. 243--246]{Griffiths-Harris}), for example.

\begin{proof}[Proof of Proposition \ref{Delta}]
Since $\GN$ has genus zero, we have an isomorphism 
$$
j_{\Gamma}:\XN \tilde{\longrightarrow}\mathbb{P}^1(\mathbb{C}),
$$
so we may view $j_{\Gamma}$ as a local coordinate for the Riemann sphere. Define the divisor $D$ by
\begin{equation}\label{divisorD}
D=\mathrm{div}((dj_{\Gamma})^{k})+(1-d+k)P_{\infty}+\sum_{\substack{\mathrm{cusps }\,P\\P\neq\infty}} kP+
\sum_{\substack{\mathrm{elliptic }\,P}} [k(1-1/e_P)]P,
\end{equation} where $P_{\infty}$ is the point ``infinity''and $[k(1-1/e_P)]$ 
is the integer part of $k(1-1/e_P),$ and if $\pi(Q)=P$ for some elliptic fixed point $Q,$ then $e_P=e_Q.$ 

Let $\omega=h(dj_{\Gamma})^k$ be a meromorphic differential $k$-form on $\Gamma \backslash \UH$.  First we claim that $\w$ corresponds a holomorphic modular form satisfying the conditions of the proposition if and only if $h\in L(D).$ Let $h\in L(D),$ then we have $\mathrm{div}(h)+D\geq 0,$ and hence,
{\allowdisplaybreaks\begin{eqnarray*}
\text{ord}_{P}(\w)+k(1-1/e_P)&\geq 0&\text{at an elliptic point,}\\
\text{ord}_{P}(\w)+k&\geq 0&\text{at a cusp,}\\
\text{ord}_{P}(\w)&\geq 0&\text{at the remaining points.}
\end{eqnarray*}}It is known \cite[Lemma 4.11]{Milne} that if $f$ is the modular form  corresponding
 to $\omega$ and $\pi(Q)=P$ under the canonical projection $\pi:\UH^* \to \Gamma \backslash \UH^*$, then the following are true:
{\allowdisplaybreaks\begin{eqnarray*}
\text{ord}_Q(f)&=&e_Q\text{ord}_{P}(\w)+k(e_Q-1)\text{ if $Q$ is an elliptic point,}\\
\text{ord}_Q(f)&=&\text{ord}_{P}(\w)+k\text{ if $Q$ is a cusp,}\\
\text{ord}_Q(f)&=&\text{ord}_{P}(\w)\text{ at the remaining points.}
\end{eqnarray*}}Therefore $\text{ord}_Q(f)\geq 0$ for all $Q\in \UH^{*},$ 
and so $f$ is holomorphic. Clearly, at $\infty$, we have $\text{ord}_{P}(\w)=d-1-k,$
 and hence $\text{ord}_{Q}(f)=d-1.$ It is also known \cite[p. 46]{Shimura} that
\begin{equation}\label{ord}
\sum \big(\text{ord}(f)/e_Q-k(1-1/e_Q)\big) =-2k+\nu_{\infty}k
\end{equation}and 
\begin{equation}\label{dim}
d=1-2k+\nu_{\infty}k+\sum [k(1-1/e_Q)],
\end{equation}where $\nu_{\infty}$ is the number of inequivalent cusps. Since $f$ has a 
zero of order $d-1$ at infinity, by comparing \eqref{ord} with \eqref{dim}, we can 
conclude that $f$ has a zero of order $e_Q\{k(1-1/e_Q)\}$ at elliptic fixed point $Q.$

Second, we want to show that $l(D)=1.$ Clearly
\begin{equation}\label{degreeD}
\mathrm{deg}(D)=-2k+(1-d+k)+(\nu_{\infty}-1)k +\sum [k(1-1/e_P)].
\end{equation}By \eqref{dim}, we then obtain that $\deg(D)=0.$ Hence by the Riemann-Roch theorem, we derive that $l(D)=1.$ Therefore there is a unique $\DN\in\Mk$ satisfying all conditions. 

We are left with proving that if $\GN$ is good, then $\Upsilon$ has real Fourier coefficients.   Notice that, in view of what we have proved above combined with Lemma \ref{trivialzero}, $\Upsilon$ is uniquely characterized by the facts that
\begin{enumerate}
\item The leading term in the Fourier expansion of $\Upsilon$ at $\infty$ is $1$, and
\item Any $f \in \Mk$ that is not a scalar multiple of $\Upsilon$ satisfies
$\mathrm{ord}_{\infty}(f) < \mathrm{ord}_{\infty}(\Upsilon)$.
\end{enumerate}
These  conditions can be rephrased as saying that among the elements of $\Mk$ with leading term $1$ in their Fourier expansion at $\infty$, $\Upsilon$ has the maximum number of zero at $\infty$ and hence its Fourier expansion at $\infty$ is $q^{d-1}+O(q^d).$  Since $\Mk$ has a basis of forms with real Fourier coefficients, it follows from simple linear algebra that $\Upsilon$ is a real linear combination of forms in this basis, and hence itself has real Fourier coefficients.
\end{proof}

\section{Proof of Theorem \ref{main}} \label{proof1}

We begin this section by recalling the Petersson inner product. If $f$ and $g$ are two modular forms in $\Mk,$ then their 
{\it Petersson inner product} is defined by 
$$ \langle f,g \rangle=\iint_{\F}f(z)\overline{g(z)}y^{2k-2}dxdy,
$$
where $\F$ denotes a fundamental domain for the action of $\GN,$ and where $z=x+iy.$ In Section \ref{sec2}, we showed that there is a divisor polynomial $P(f,\JN)$ (resp. $P(g,\JN)$) attached to $f$ (resp. $g$) such that $f=P(f,\JN)\DN$ (resp. $g=P(g,\JN)\DN$). So we derive that
\begin{equation}\label{Petersson}
\langle f,g \rangle=\iint_{\F}P(f,\JN(z))\overline{P(g,\JN(z))}|\DN(z)|^2y^{2k-2}dxdy.
\end{equation}It is well-known that if $f(z)g(z)$ vanishes at every cusp of $\GN,$ then the integral \eqref{Petersson} converges.

We denote by $\overline{\F}$ the closure of the acceptable fundamental domain $\F$ for $\GN.$
Write 
\begin{eqnarray} \label{Ai}
A=\amalg_{i=1}^n\{A_i\}
\end{eqnarray}
 where the $A_i$ are smooth (rectifiable) paths (i.e. smooth one-dimensional real manifolds which may be
non-compact or with boundary) such that $j_{\GN}|_{\overline{A}_i}$ is injective.  We then write $\overline{\F}=\bigcup_{i=1}^{n+2}\F_i$ where the $\F_i$ are the closures of domains in $\mathfrak{H}$ of positive area with respect to the hyperbolic metric that satisfy the following properties:
\begin{enumerate}
\item{If $i \neq j,$ $\F_i \cap \F_j$ is a smooth curve that is rectifiable with respect to the hyperbolic metric.}
\item{For $1 \leq i \leq n$ we have $A_i \subset \F_i,$ and $A_i \cap \F_j$ is either empty or a single point for $i \neq j.$}
\item{$\partial \F_i - \overline{A_i}$ is contained in the interior of $\F.$}
\item{$\F_{n+1}= \overline{\F} \cap \{z: \Re(z) \leq 0\} \cap \{z: \Im(z)\geq y_0\}$ and $\F_{n+2}=\overline{\F} \cap \{z: \Re(z) \geq 0\} \cap \{z: \Im(z)\geq y_0\}$ where $y_0$ is defined as \eqref{y0}.}
\item{$j_{\GN}|_{\F_i}$ is injective.}
\end{enumerate}
One can always write $\overline{\F}$ as such a union; we do this to avoid pathologies with certain line integrals which we will now introduce 

For $z=x+iy\in\overline{\F}_i,$ define a path $\gamma_i(z)$ by
\begin{equation}\label{path}
\gamma_i(z):=\{\JN^{-1}(\sigma+i\Im(\JN(z)))\,:\,\sigma\in(\Re(\JN(z)),\infty)\} 
\end{equation}traversed from $\infty$ to $z.$ Let $\Phi_i(z)$ be the path integral
\begin{equation}\label{Phi}
\Phi_i(z):=\int_{\gamma_i(z)}\Re\Big(P(f,\JN(s))\overline{P(g,\JN(s))}|\DN(s)|^{2}t^{2k-2}\Big)dt,
\end{equation}where $\Im(s)=t.$ Then we obtain the following result.

\begin{proposition}\label{integral}
Let $\GN$ be good for the weight $k$ of genus zero and $\F$ be the acceptable fundamental domain for $\GN.$ Let $f$ and $g$ be in $\Mk$ such that $fg$ vanishes at every cusp of $\GN.$ Then we have that 
\begin{equation}\label{lowerarcs}
\iint_{\F}\Re\Big(P(f, \JN(z))\overline{P(g,\JN(z))}|\DN(z)|^{2}y^{2k-2}\Big)dxdy=-\sum_{i=1}^n\int_{A_i}\Phi_i(z)dx,
\end{equation} where $z:=x+iy$, and the $A_i$ are the smooth paths in (\ref{Ai})  traversed with counterclockwise orientation.
\end{proposition}

We will prove Proposition \ref{integral} after the proof of Theorem \ref{main}:

\begin{proof}[Proof of Theorem \ref{main}.]
 Let $\frak{a}_1, \dots, \frak{a}_m$ be the zeros of $P(\Ek,\JN)$ along $\F \cap \JN^{-1}([\JN(-\frac{h}{2}+iy_0), \infty))$ that have odd multiplicity and let $\frak{b}_1, \dots,\frak{b}_{c(\GN, \F)}$ be the points along the lower arc $\F\cap A$ where $\frac{\partial y}{\partial \JN}$ changes sign. If $m+c(\F,\Gamma) \geq \mathrm{deg} (P(\Ek,X))$, we are done, so we assume for a contradiction that $m+c(\GN,\F) <\mathrm{deg}(P(\Ek,X))$.  Let
\begin{equation}\label{Q}
Q(X):=\prod_{i=1}^m(X-\JN(\frak{a}_i))\prod_{j=1}^{c(\Gamma,\mathcal{F})}(X-\JN(\frak{b}_j)).
\end{equation} 
Then $Q(\JN)$ is a polynomial in $\JN$ with real coefficients.

Note that $Q(\JN)\DN \in M_{2k}^{\infty}(\Gamma)$.  In fact, by Lemma \ref{Delta} and the construction of $Q(X)$, we have that $Q(\JN)\DN \in M_{2k}(\Gamma)$ and $Q(\JN)\DN$ vanishes at $\infty$.  Since any modular form in $\Mk$ that vanishes at infinity is orthogonal to $\Ek$ with respect to the Petersson inner product,  we have
\begin{equation}\label{innerzero}
0=\langle \Ek, Q(\JN)\DN\rangle.
\end{equation}

On the other hand, we have that for $z=x+iy,$
{\allowdisplaybreaks\begin{align*}
\langle \Ek, Q(\JN)\DN\rangle=&\iint_{\F}\Re\Big(P(\Ek,\JN(z))\overline{Q(\JN(z))}|\DN(z)|^{2}\Big)y^{2k-2}dxdy\\&
+i\iint_{\F}\Im\Big(P(\Ek,\JN(z))\overline{Q(\JN(z))}|\DN(z)|^{2}\Big)y^{2k-2}dxdy.
\end{align*}}By Proposition \ref{integral} we know that
\begin{align}
\iint_{\F}&\Re\Big(P(\Ek, \JN(z))\overline{Q(\JN(z))}|\DN(z)|^{2}y^{2k-2}\Big)dxdy\nonumber\\
=&-\sum_{i=1}^n\int_{A_i}\Big(\int_{\gamma_i(z)}P(\Ek,\JN)Q(\JN)(|\DN|^2y^{2k})(\JN)y^{-2}(\JN)y'(\JN)d\JN\Big)dx.\label{maineqn}
\end{align}
where $\gamma_i(z)$ is as in (\ref{path}).  Here we have used the fact that $\JN|_{\mathcal{F}_i}$ is injective, so we may view $y$ as a function of $\JN$.  Our path $\gamma_i(z)$ was chosen so that whenever $z \in A$, the function $\JN(s)$ is real for all $s \in \gamma_i(z)$.  
Moreover, by construction of $Q(\JN)$, we have that for each $z \in A_i$,  
$$
P(\Ek,\JN(s))Q(\JN(s))y'(\JN(s))
$$
is nonnegative or nonpositive for $s \in \gamma_i(z)$, and its sign does not depend on $i$.  
It follows that the integral \eqref{maineqn} can not be zero. This is the desired contradiction.
\end{proof}

We now give a proof of Proposition \ref{integral}:

\begin{proof}[Proof of Proposition \ref{integral}]
Let $ \delta\in \{1,2\}.$ We define 
$$
\F_{n+\delta}^T:=\F_{n+\delta} \cap \{z: \Im(z) \leq T\}
$$ for $T \gg 0,$ and let $L_{n+\delta}^T$ denote the horizontal lines 
$$
L_{n+\delta}^T=\F_{n+\delta}^T\cap\{z\,:\,\Im(z)=T\}.
$$ 
Note that 
\begin{eqnarray} \label{dif}
\frac{\partial \Phi_i(z)}{\partial y}=\Re\Big(P(f, \JN(z))\overline{P(g, \JN(z))}|\DN(z)|^{2}y^{2k-2}\Big).
\end{eqnarray}
We apply Green's theorem to each $\mathcal{F}_i$ separately, use (\ref{dif}), and then sum the contributions to obtain  
{\allowdisplaybreaks\begin{align} \label{fonow}
\iint_{\F}&\Re\Big(P(f, \JN(z))\overline{P(g, \JN(z))}|\DN(z)|^{2}y^{2k-2}\Big)dxdy\\ \nonumber
=&-\sum_{\delta=1}^{2}\lim_{T\rightarrow \infty}\int_{L_{n+\delta}^T}\Phi_{n+\delta}(z)dx-\sum_{i=1}^n\int_{A_i}\Phi_i(z)dx.
\end{align}}

\noindent Here the limit exists because $fg$ vanishes at $\infty$ (and hence has rapid decay).  The contribution of the vertical lines contained in $\left\{z: \Re(z)=\pm \frac{h}{2}\right\} \cap \overline{\F}$ is zero because $dx=0$ there. We now claim that
\begin{equation}\label{Phi-T}
\lim_{T\rightarrow \infty}\int_{L_{n+\delta}^T}\Phi_{n+\delta}(z)dx=0;
\end{equation}
it is easy to see that \eqref{Phi-T} in conjunction with \eqref{fonow} will finish the proof of the proposition.

In view of the fact that $\int_{L_{n+\delta}^T}dx=h/2<\infty$ for $T \gg 0$, in order to prove \eqref{Phi-T}, it is enough to show that 
\begin{eqnarray} \label{lastpart}
\Phi_{n+\delta}(z)=\int_{\gamma_{n+\delta}(z)}\Re\Big(P(f, \JN(s))\overline{P(g, \JN(s))}|\DN(s)|^{2}t^{2k-2}\Big)dt \longrightarrow 0
\end{eqnarray}as $y\rightarrow\infty$ while $z \in \F_{n+\delta}$. Here $t=\Im(s).$ Because $fg$ vanishes at every cusp, we have that  
$$\Re\Big(P(f, \JN(s))\overline{P(g, \JN(s))}|\DN(s)|^{2}t^{2k}\Big) =\Re\left(f(s)\overline{g(s)}t^{2k}\right)<C $$ 
for every $s \in \mathcal{F}$, where $C>0$ is some constant not depending on $s$.  Therefore, to prove \eqref{lastpart} it suffices to show that
 \begin{eqnarray} \label{pathzero}
 \int_{\gamma_{n+\delta}(z)}t^{-2}dt \longrightarrow 0
 \end{eqnarray}
 as $y \to \infty$ while $z \in \F_{n+\delta}$.
Note that since $\JN(x+iy)\sim e^{-2\pi i(x+iy)}$ as $y \to \infty$, we have 
$$
\JN^{-1}(x+iy)\sim \frac{1}{2\pi}\left(i\log \sqrt{x^{2}+y^{2}}-\text{Arg}(x+iy)\right).
$$ 
where $\log$ denotes the principal branch of the natural logarithm.  
Thus, as $y \to \infty$ while $z \in \F_{n+\delta}$,
 {\allowdisplaybreaks\begin{align*}
\int_{\gamma_{n+\delta}(z)}t^{-2}dt
&\sim \int_{\infty}^{\Re(\JN(z))}\log \sqrt{x^{2}+\Im(\JN(z))^{2}}^{-2}(\log \sqrt{x^{2}+\Im(\JN(z))^{2}})'dx\\
&\sim\int_{\infty}^{\Re(\JN(z))^2+\Im(\JN(z))^2}\frac{dx}{x\log^{2}x}
\\&\sim \frac{1}{\log(\Re(\JN(z))^2+\Im(\JN(z))^2)}\to 0,
\end{align*}}as $z \to \infty,$ since$\JN(z) \sim e^{-2 \pi i(x+iy)}.$ This implies \eqref{lastpart}, and in turn the proposition.
\end{proof}




\begin{thebibliography}{99}

\bibitem{AAR}
G. E. Andrews, R. Askey, and R. Roy, \emph{Special functions}, Encyclopedia of Mathematics \textbf{71}, Cambridge Univ. Press, 
Cambridge, 1999.

\bibitem{Boylan}
M. Boylan, \emph{Swinnerton-Dyer type congruences for certain Eisenstein series}, Contemporary Mathematics \textbf{291} (2001), 93--108.


\bibitem{Diamond-Im}
F. Diamond and J. Im, \emph{Modular forms and modular curves}, Canadian Math. Soc. Conference Proceedings \textbf{17} (1995), 39--133.


\bibitem{Getz}
J. Getz, \emph{A generalization of a theorem of Rankin and Swinnerton-Dyer on zeros of modular forms}, Proc. American Math. Soc. \textbf{132} (2004), 2221--2231.

\bibitem{Griffiths-Harris}
P. Griffiths and J. Harris, \emph{Principles of algebraic geometry}, A Wiley-Interscience publication, New York, 1978


\bibitem{Milne}
J. S. Milne, \emph{Modular functions and modular forms}, lecture notes, University of Michigan, 1997.

\bibitem{Ono}
K. Ono, \emph{The web of modularity: Arithmetic of the coefficients of modular forms and $q$-series}, CBMS, \textbf{102}, American Math. Soc., Providence, Rhode Island, 2004. 

\bibitem{Ono-Bring}
K. Ono and K. Bringmann, \emph{Identities for traces of singular moduli}, Acta Arith. \textbf{119} (2005), 317--327.

\bibitem{Rankin-Swinnerton}
F. K. C. Rankin and H. P. F. Swinnerton-Dyer, \emph{On the zeros of Eisenstein series}, Bull. London Math. Soc. \textbf{2} (1970),  169--170.


\bibitem{Rankin-Poincare}
R. A. Rankin, \emph{The zeros of certain Poincar\'e series}, Compositio Math. \textbf{46} (1982), 255--272.

\bibitem{Rudnick}
Z. Rudnick, \emph{On the asymptotic distribution of zeros of modular forms}, IMRN, No. 34 (2005), 2059--2076

\bibitem{Schoeneberg}
B. Schoeneberg, \emph{Elliptic modular functions}, Springer--Verlag, New York, Heidelberg, Berlin, 1970.

\bibitem{Shimura}
G. Shimura, \emph{Introduction to the arithmetic theory of automorphic functions}, Princeton Univ. Press, Princeton, 1971.

\bibitem{Verrill}
H. A. Verrill, \emph{Fundamental domain drawer, Java}, http://www.math.lsu.edu/verrill/fundomain

\end{thebibliography}
\end{document}